\documentclass[11pt]{article}%
\usepackage{amssymb, amsmath, latexsym, mathrsfs, verbatim, calc}
\usepackage[latin1]{inputenc}
\usepackage{amsmath}
\usepackage{amsfonts}
\usepackage{amssymb}
\usepackage{graphicx}
\usepackage{color}%
\providecommand{\U}[1]{\protect\rule{.1in}{.1in}}
\def \bop {\noindent\textbf{Proof. }}
\def \eop {\hbox{}\nobreak\hfill
\vrule width 2mm height 2mm depth 0mm
\par \goodbreak \smallskip}

\newtheorem{theorem}{Theorem}

\newtheorem{definition}[theorem]{Definition}

\newtheorem{lemma}[theorem]{Lemma}

\newtheorem{proposition}[theorem]{Proposition}
\newtheorem{remark}[theorem]{Remark}

\newcommand{\R}{\mathbb{R}}

\def\R{\mathbb{R}}

\def\U{\mathbb{U}}
\begin{document}

\title{Existence of an Optimal Control for a coupled FBSDE with a non degenerate
diffusion coefficient}
\author{K. BAHLALI\thanks{Université de Toulon, IMATH, EA $2134$, $83957$ La Garde
cedex, France.}
\and O. KEBIRI\thanks{Tlemcen University, Laboratory of Probability and Statistics,
Algeria.}
\and B. MEZERDI\thanks{Université Mohamed Khider, Biskra, Algerie.}
\and A. MTIRAOUI \footnotemark[1] }
\date{}
\maketitle

\vskip -0.5cm


\vskip -0.4cm

\noindent\textbf{Abstract}: We a controlled system  driven by a coupled
forward-backward stochastic differential equation (FBSDE) with a non
degenerate diffusion matrix. The cost functional is defined by the solution of
the controlled backward stochastic differential equation (BSDE), at the
initial time. Our goal is to find an optimal control which minimizes the cost
functional. The method consists to construct a sequence of approximating
controlled systems for which we show the existence of a sequence of feedback
optimal controls. By passing to the limit, we establish the existence of a
relaxed optimal control to the initial problem. The existence of a strict
control follows from the Filippov convexity condition. Our result improve in
some sense those of \cite{bgm,BLRT}.

\noindent\textbf{Keywords:} Optimal control, forward-backward stochastic
differential equations, stochastic control, Hamilton-Jacobi-Bellman equation,
relaxed control, strict control.


\section*{Introduction}

Stochastic optimal control has interested many researchers, both for its
theoretical aspect as well as for its applications in real world problems.
There is a vast literature dealing with optimal control for systems driven by
stochastic differential equations (SDE) and/or forward-backward SDEs (FBSDE),
and various aspects were studied. The principal developments concern the
existence of optimal control, Pontryagin's maximum principle (or necessary
optimality conditions) and Bellman's principle (also called dynamic
programming principle), see e.g. \cite{bgm, bi, bl, BLRT, enj, fs, hl,
Krylov2, Li, LW, PW}. Closer to or concern here, the existence of an optimal
control for a system driven by an FBSDE was established in \cite{bgm} and
\cite{BLRT} by different methods. In \cite{bgm}, the approach consists to
directly show the existence of a relaxed control by using a compactness method
and the Jakubowsky $S$-topology. In \cite{BLRT} the authors pass by the
Hamilton Jacobi Bellman (HJB) equation associated with the control problem.
This allows them to construct a sequence of optimal feedback controls. Then,
they pass to the limit and use the result of \cite{enj} in order to get the
existence of a relaxed optimal control. In both papers \cite{bgm} and
\cite{BLRT} the Filippov convexity condition is used in order to get the
existence of a strict optimal control. It should be noted that in \cite{bgm}
and \cite{BLRT} the controlled system is driven by a decoupled FBSDE.

The aim of the present paper is to extend the results of \cite{bgm, BLRT}, to
a coupled FBSDE. To begin, let us give a description of our problem.

Let $T>0$ be a finite horizon, $t\in\lbrack0,\ T]$ and $(\Omega,\,\mathcal{F}%
,\,\mathbb{P},\,(\mathcal{F}_{t}))$ be a filtered probability space satisfying
the usual conditions. Let $W$ be a $d$-dimensional Brownian motion with
respect to the filtration $(\mathcal{F}_{t})$. Let $\mathbb{U}$ be a compact
metric space. We define the deterministic functions $b$, $\sigma$, $f$ and
$\Phi$ by
\[
b:\mathbb{R}^{d}\times\mathbb{R}\times\mathbb{U}\longmapsto\mathbb{R}^{d},
\]%
\[
\sigma:\mathbb{R}^{d}\times\mathbb{R}\longmapsto\mathbb{R}^{d\times d},
\]%
\[
f:\mathbb{R}^{d}\times\mathbb{R}\times\mathbb{R}^{d}\times\mathbb{U}%
\longmapsto\mathbb{R},
\]%
\[
\Phi:\mathbb{R}^{d}\longmapsto\mathbb{R}.
\]

We consider the following controlled coupled FBDSE defined for $s\in\lbrack
t,\,T]$ by:%

\begin{equation}
\left\{
\begin{array}
[c]{l}%
dX_{s}^{t,x,u}=b(X_{s}^{t,x,u},Y_{s}^{t,x,u},u_{s})ds+\sigma(X_{s}%
^{t,x,u},Y_{s}^{t,x,u})dW_{s},\\
X_{t}^{t,x,u}=x,\,\\
dY_{s}^{t,x,u}=-f(X_{s}^{t,x,u},Y_{s}^{t,x,u},Z_{s}^{t,x,u},u_{s}%
)ds+Z_{s}^{t,x,u}dW_{s}+dM_{s}^{t,x,u},\\
M_{t}^{t,x,u}=0,\ \ \langle M^{t,x,u},\ W\rangle_{s}=0,\ \ Y_{T}^{t,x,u}%
=\Phi(X_{T}^{t,x,u})
\end{array}
\right.  \label{equa}%
\end{equation}
where, $X^{t,x,u},\,Y^{t,x,u},\,Z^{t,x,u}$ are $(\mathcal{F}_{t})$-adapted
square integrable processes and $M^{t,x,u}$ is an $(\mathcal{F}_{t})$-adapted
square integrable martingale, which is orthogonal to $W.$ The control variable
$u$ is an $\mathcal{F}_{t}-$ adapted process with values in $\mathbb{U}$. It
should be noted that the filtered probability space and the Brownian motion
may change with the control $u$.

On $\nu:=(\Omega,\mathcal{F} ,\mathbb{P} ,\mathcal{F}_{t} ,W )$, we define the
following spaces of processes, for $m\in\mathbb{N}^{*}$ and $t\in[0,T)$,

\begin{itemize}
\item $\mathcal{S}^{2}_{\nu}(t,T;\mathbb{R}^{m})$ \ denote the set of
$\mathbb{R}^{m}$-valued, $\mathcal{F}_{t} $-adapted, continuous processes
$(X_{s},s\in[t,T])$ which satisfy $\mathbb{E}(\sup_{t\leq s\leq T}|X_{s}
|^{2})<\infty$,

\item $\mathcal{H}^{2}_{\nu}(t,T;\mathbb{R}^{m})$ \ is the set of
$\mathbb{R}^{m}$-valued, $\mathcal{F}_{t} $-predictable processes $(Z
_{s},s\in[t,T])$ which satisfy $\mathbb{E}[\int_{t}^{T}|Z_{s}|^{2}ds]<\infty$,

\item $\mathcal{M}_{\nu}^{2}(t,T;\mathbb{R}^{m})$ \ denotes the set of all
$\mathbb{R}^{m}$-valued, square integrable càdlàg martingales $M=(M_{s}%
)_{s\in\lbrack t,T]}$ with respect to $\mathcal{F}_{t}$, with $M_{t}=0$.
\end{itemize}

\begin{definition}
\label{admicont} A solution of the FBSDE (\ref{equa}) is a process
$(X^{t,x,u},Y^{t,x,u},Z^{t,x,u},M^{t,x,u})\in\mathcal{S}_{\nu}^{2}%
(t,T;\mathbb{R}^{d})\times\mathcal{S}_{\nu}^{2}(t,T;\mathbb{R})\times
\mathcal{H}_{\nu}^{2}(t,T;\mathbb{R}^{d})\times\mathcal{M}_{\nu}%
^{2}(t,T;\mathbb{R}^{d})$ satisfying equation (\ref{equa}).
\end{definition}

\begin{definition}
1) A strict control is an $\mathcal{F}_{t}$-progressively measurable processes
$(u_{s},s\in\lbrack t,T])$ with values in $\mathbb{U}$, such that FBSDE
(\ref{equa}) has a solution in \ $\mathcal{S}_{\nu}^{2}(t,T;\mathbb{R}%
^{d})\times\mathcal{S}_{\nu}^{2}(t,T;\mathbb{R})\times\mathcal{H}_{\nu}%
^{2}(t,T;\mathbb{R}^{d})\times\mathcal{M}_{\nu}^{2}(t,T;\mathbb{R}^{d})$. We
denote $\mathcal{U}_{\nu}(t)$ the set of all strict controls.

2) A relaxed control is an $\mathcal{F}_{t}$-progressively measurable
processes $(\mu_{s},s\in\lbrack t,T])$ with values in the space $\mathcal{P}%
(\mathbb{U)}$ of probability measures in $\mathbb{U}$, such that FBSDE
(\ref{equa}) has a solution in \ $\mathcal{S}_{\nu}^{2}(t,T;\mathbb{R}%
^{d})\times\mathcal{S}_{\nu}^{2}(t,T;\mathbb{R})\times\mathcal{H}_{\nu}%
^{2}(t,T;\mathbb{R}^{d})\times\mathcal{M}_{\nu}^{2}(t,T;\mathbb{R}^{d})$. We
denote $\mathcal{R}_{\nu}(t)$ the set of all relaxed controls.
\end{definition}

The cost functional, which will be minimized, is defined for $u \in
\mathcal{U}_{\nu}(t)$ \ by:
\begin{equation}
\label{j}J(t,x,u):=Y^{t,x,u}_{t}.
\end{equation}

An $\mathcal{F}_{t}$-adapted control $\widehat{u}$ is called optimal if it
minimizes $J$, that is:%

\[
Y^{t,x,\widehat{u}}_{t} = \mbox{essinf}\left\{  Y^{t,x,u}_{t} ,\,
\ u\in\mathcal{U}_{\nu}(t)\right\}  .
\]

If $\widehat{u}$ belongs to $\mathcal{U}_{\nu}(t)$, we then say that
$\widehat{u}$ is an optimal strict control.

The value function $V$ is defined by:
\begin{equation}
\label{original}V(t,x):=Y^{t,x,\hat{u}}_{t}= \mbox{essinf}\left\{  J(t,x,u)
,\, \ u\in\mathcal{U}_{\nu}(t)\right\}  .
\end{equation}

Our objective is to establish the existence of a strict optimal control for
the problem (\ref{equa})--(\ref{original}). To this end, we follow the method
developed in \cite{BLRT}: we approximate the controlled FBSDE (\ref{equa}) by
a sequence of FBSDEs with smooth data $b_{\delta},\sigma_{\delta},f_{\delta}$
and $\Phi_{\delta}$ and consider the new value function $V^{\delta}$ which is
associated to the FBSDE, with these smooth data. According to Krylov
\cite{Krylov1} (Theorems 6.4.3 and 6.4.4), $V^{\delta}$ is sufficiently smooth
and satisfies a Hamilton-Jacobi-Bellman equation. Since all admissible
controls take their values in a compact set, we then deduce the existence of a
feedback optimal control $u^{\delta}$. Next, we prove that the sequence
$V^{\delta}$ converges uniformly to a function $V$, which is the value
function of our initial control problem. Comparing with \cite{BLRT}, The first
difficulty is related to the fact that: if we consider the usual definition of
admissible controls, then the uniform Lipschitz condition on the coefficients
[assumption (A1)] is not sufficient  to prove the existence of a unique
solution to equation (\ref{equa}) for an arbitrary duration. This fact is well
explained in \cite{antonelli} and two examples are given. For this reason, we
had to change the usual definition of admissible controls and adopt the
definition \ref{admicont} above. In order to ensure that our definition has a
sense, that is the set of admissible controls is not empty, we moreover assume
throughout this paper that the diffusion matrix $\sigma$ is non degenerate
[assumption (A3)]. In this case, the set of admissible controls contains the
constants. Indeed, if the control $u$ is constant then according to \cite{del}
the FBSDE (\ref{equa}) has a unique solution in \ $\mathcal{S}_{\nu}%
^{2}(t,T;\mathbb{R}^{d})\times\mathcal{S}_{\nu}^{2}(t,T;\mathbb{R}%
)\times\mathcal{H}_{\nu}^{2} (t,T;\mathbb{R}^{d})\times\mathcal{M}_{\nu}%
^{2}(t,T;\mathbb{R}^{d})$. The second difficulty concerns the uniform estimate
of the variable $Z$ and the stability of solutions. The later require, in our
situation, a harder computation which combines PDEs techniques and FBSDEs
arguments. Assuming that the diffusion matrix is non degenerate and the
coefficients are uniformly Lipschitz in $(x, y, z)$ and continuous in the
control $u$, we establish the existence of an admissible feedback control $u$
by decoupling the FBSDE (\ref{equa}) and by using the results of
\cite{bahlali} and \cite{PP1}. The method consists to construct a sequence of
approximating controlled systems for which we show the existence of a sequence
of feedback optimal controls. By passing to the limit, we establish the
existence of a relaxed optimal control to our initial problem. The existence
of a strict control follows from the Filippov convexity condition. Note that
when the control enters the diffusion coefficient $\sigma$, we obtain a SDE
with a measurable diffusion matrix and, in this case, the uniqueness of
solution (even in law sense) may fails. Indeed, we know from \cite{krylovspa}
that when the diffusion coefficient is merely measurable, then even the
uniqueness in law fails in general for Itô's forward SDE in dimension strictly
greater than 2, see \cite{krylovspa} for more details. This explains why we
consider only the case when the control does not enter the diffusion coefficient.

The paper is organized as follows. In section 1, we introduce some notations,
the controlled system, and the assumptions. In section 2, we present the cost
functional and the value function which satisfies the Hamilton-Jacobi-Bellman
equation. In section 3, we give the main result and its proof. This section
contains two subsections. The first one is devoted to the study of the
approximating control problem, together with its associated HJB equation. In
the second subsection, we prove our main result.


\vskip 0.2cm\noindent\textbf{Assumption (A)}.

\begin{itemize}
\item \textbf{(A1)} \textbf{1)} \ There exists $K>0$ such that for any
$u\in\mathbb{U}$, $(x,y,z)\mbox{ and }(x^{\prime},y^{\prime},z^{\prime}%
)\in\mathbb{R}^{d}\times\mathbb{R}\times\mathbb{R}^{d}$
\begin{align*}
&  \left\vert \sigma(x,y)-\sigma(x^{\prime},y^{\prime})\right\vert \leq
K(\left\vert x-x^{\prime}\right\vert +\left\vert y-y^{\prime}\right\vert ),\\
&  |\Phi(x)-\Phi(x^{\prime})|\leq K|x-x^{\prime}|,\\
&  |b(x,y,u)-b(x^{\prime},y^{\prime},u)|\leq K(|x-x^{\prime}|+|y-y^{\prime
}|),\\
&  |f(x,y,z,u)-f(x^{\prime},y^{\prime},z^{\prime},u)|\leq K(|x-x^{\prime
}|+|y-y^{\prime}|+|z-z^{\prime}|).
\end{align*}
\textbf{2)} \ The functions $b$, $\sigma$, $f$ and $\Phi$ are bounded.

\item \textbf{(A2)} \ For every $(x,y,z)\in\mathbb{R}^{d}\times\mathbb{R}%
\times\mathbb{R}^{d},$ the functions $b(x,y,.)$ and $f(x,y,z,.)$ are
continuous in $u$.

\item \textbf{(A3)} There exists $\lambda>0$ such that for every
$(t,x,y)\in\lbrack0,T]\times\mathbb{R}^{d}\times\mathbb{R}$,
\[
\forall\zeta\in\mathbb{R}^{d}~~~~\langle\zeta,\sigma(t,x,y)\zeta\rangle
\geq\lambda|\zeta|^{2}%
\]

\end{itemize}

When the control $u$ is constant, one can show (as in \cite{del}) that under
assumptions (A1) and (A3), equation (\ref{equa}) has a unique solution
\ $(X^{t,x,u},Y^{t,x,u},Z^{t,x,u},M^{t,x,u})$ in the space \ $\mathcal{S}%
_{\nu}^{2}(t,T;\mathbb{R}^{d})\times\mathcal{S}_{\nu}^{2}(t,T;\mathbb{R}%
)\times\mathcal{H}_{\nu}^{2}(t,T;\mathbb{R}^{d})\times\mathcal{M}_{\nu}%
^{2}(t,T;\mathbb{R}^{d}).$

\vskip 0.2cm The following assumption \textbf{{(H)}} will be called the
convexity assumption. \label{convexity}
\[
\mathbf{{(H)}} \qquad\left\{
\begin{array}
[c]{l}%
\mbox{For all $(x,y)\in\R^d\times\R$ the following set is convex:}\\
\{ ((\sigma\sigma^{*})(x,y),w(\sigma\sigma^{*})(x,y),b(x,y,u),f(x,y,w\sigma
(x,y),u))\\
|(u,w)\in\mathbb{U} \times\bar{B}_{C}(0)\}\; ,
\end{array}
\right.
\]
where $\bar{B}_{C}(0)\ \subset\mathbb{R}^{d}$ is the closed ball around 0 with
radius $C$.

\vskip0.2cm The following lemma can be proved as Lemma 4 of \cite{BLRT}. For
completeness, we give its proof in the appendix.

\begin{lemma}
\label{H1} For $(x,y,w,\theta,u)\in\mathbb{R}^{d}\times\mathbb{R}%
\times\mathbb{R}^{d}\times\mathbb{R}\times\mathbb{U}$, set
\[
\Sigma(x,y,w,\theta)=\left(
\begin{array}
[c]{ll}%
\sigma(x,y) & 0\\
w\sigma(x,y) & \theta
\end{array}
\right)  \;\;\mbox{ and }\;\;\beta(x,y,w,u)=\left(
\begin{array}
[c]{c}%
b(x,y,u)\\
-f(x,y,w\sigma(x,y),u)
\end{array}
\right)  \;.
\]
Under assumption \textbf{(H)} we have
\[
\overline{co}\{((\Sigma\Sigma^{\ast})(x,y,w,0),\beta(x,y,w,u))|(u,w)\in
\newline U\times\bar{B}_{C}(0)\}
\]%
\[
\subset\{((\Sigma\Sigma^{\ast})(x,y,w,\theta),\beta(x,y,w,u)|(u,w,\theta
)\in\mathbb{U}\times\bar{B}_{C}(0)\times\lbrack0,K]\}
\]
where, for any set $\mathbb{E}$, $co(E)$ denotes the convex hull of $E$.
\end{lemma}

\noindent\textit{The Hamilton-Jacobi-Bellman equation}

Let $\mathbb{S}^{d}$ denotes the space of the symmetric matrices in
$\mathbb{R}^{d^{2}}$. For a function $V$, we denote by $\nabla_{x}V$ the
gradient and $\nabla_{xx}V$ the Hessian of the matrix $V$. Let $H$ be the real
function defined on $\mathbb{R}^{d}\times\mathbb{R}\times\mathbb{R}^{d}%
\times\mathbb{S}^{d}\times\mathbb{U}$ by:
\begin{equation}
H(x,y,p,A,u):=\frac{1}{2}\mbox{tr}\left(  (\sigma\sigma^{\ast})(x,y)A\right)
+b(x,y,u)p+f(x,y,p\;\sigma(x,y),u) \label{H}%
\end{equation}
According to Li and Wei \cite{LW}, the value function $V(t,x)$, defined by
\eqref{original}, solves the following Hamilton-Jacobi-Bellman equation in the
sense of viscosity solutions.
\begin{equation}
\left\{
\begin{array}
[c]{l}%
\displaystyle\frac{\partial}{\partial t}V(t,x)+\inf_{u\in\mathbb{U}%
}H(x,V(t,x),\nabla_{x}V(t,x),\nabla_{xx}V(t,x),u)=0,\;(t,x)\in\lbrack
0,T]\times\mathbb{R}^{d},\\
V(T,x)=\Phi(x),\;x\in\mathbb{R}^{d},
\end{array}
\right.  \label{hjb}%
\end{equation}


\section{The main results}

\begin{theorem}
Assume that assumptions $\mathbf{{(A)}}$ and $\mathbf{{(H)}}$ are satisfied.
Then there exists a strict control which solves the problem \eqref{equa} and
\eqref{original} in some reference stochastic system $\bar\nu=(\bar\Omega,
\bar{\mathcal{F}} ,\bar{\mathbb{P}}, \bar{(\mathcal{F}_{t})},\bar W)$.
\end{theorem}

\noindent To prove this theorem, we approximate the controlled FBSDE
(\ref{equa}) by a sequence of FBSDEs, with smooth data $b_{\delta}%
,\sigma_{\delta},f_{\delta}$ and $\Phi_{\delta}$ and consider a new value
function $V^{\delta},$ which is associated to the FBSDE with these smooth
data. According to Krylov \cite{Krylov1} (Theorems 6.4.3 and 6.4.4),
$V^{\delta}$ is sufficiently smooth and satisfies a Hamilton-Jacobi-Bellman
equation. Since all admissible controls take their values in a compact set, we
then deduce the existence of a feedback control $u^{\delta}$. Next, we prove
that the sequence $V^{\delta}$ converges uniformly to a function $V$ which is
the value function of our initial control problem.


\section{Proof}

\subsection{Construction of an approximating Control Problem}

For an arbitrary dimension $m\geq1$ we let $\varphi:\mathbb{R}^{m}%
\mathbb{\rightarrow R}$ be a non-negative smooth function on the Euclidean
space $\mathbb{R}^{m}$ whose support is included in the unit ball of
$\mathbb{R}^{m}$ and $\int_{\mathbb{R}^{m}}\varphi\left(  \xi\right)  d\xi=1.$
Let $g:\mathbb{R}^{m}\rightarrow\mathbb{R}$ be a uniformly Lipshitz function.
We set
\[
g_{\delta}\left(  \xi\right)  =\delta^{-m}\int_{\mathbb{R}^{m}}g\left(
\xi-\xi^{^{\prime}}\right)  \varphi\left(  \delta^{-1}\xi^{^{\prime}}\right)
d\xi^{^{\prime}},\quad\xi\in\mathbb{R}^{m},\,\delta>0.
\]
The following proposition is classic and can be easily checked.

\begin{proposition}
\label{mollif} $\mbox{ for every }\xi,\xi^{\prime}\in\mathbb{R}^{m},\,
\delta,\delta^{\prime}>0,$ we have: \newline$(i)\,\left\vert g_{\delta}\left(
\xi\right)  -g\left(  \xi\right)  \right\vert \leq L_{g}\delta$\newline%
$(ii)\,\left\vert g_{\delta}\left(  \xi\right)  -g_{\delta^{\prime}}\left(
\xi\right)  \right\vert \leq L_{g}\vert\delta-\delta^{\prime}\vert,$%
\newline$(iii)\left\vert g_{\delta}\left(  \xi\right)  -g_{\delta}\left(
\xi^{\prime}\right)  \right\vert \leq L_{g}\vert\xi-\xi^{\prime}\vert,\, \,$
\newline where $L_{g}$ denotes the Lipschitz constant of $g$.
\end{proposition}

\begin{definition}
For each $\delta\in(0,1]$ we denote by $b_{\delta},\sigma_{\delta},f_{\delta}$
and $\Phi_{\delta}$ the mollifiers of the functions $b,\sigma,f$ and $\Phi,$
respectively, introduced in the second Section , with $g=b\left(  .,v\right)
,$ $\sigma\left(  .\right)  ,$ $f\left(  .,v\right)  $ and $\Phi\left(
.\right)  .$
\end{definition}


\noindent\textit{The approximating Hamilton-Jacobi-Bellman equation}

Assume that $\mathbf{{(A)}}$ is satisfied and let $\delta\in(0,1]$ be an
arbitrarily fixed number. For $\left(  x,y,p,A,v\right)  \in\mathbb{R}%
^{d}\times\mathbb{R}\times\mathbb{R}^{d}\times\mathbb{S}^{d}\times\mathbb{U}$,
we define the function $H^{\delta}$ by:
\begin{equation}
H^{\delta}\left(  x,y,p,A,v\right)  =\frac{1}{2}\left(  \mbox{tr}\left(
(\sigma_{\delta}\sigma_{\delta}^{\ast})\left(  x,y\right)  \right)  A\right)
+b_{\delta}\left(  x,y,v\right)  p+f_{\delta}\left(  x,y,p\sigma_{\delta
}\left(  x,y\right)  ,v\right)  . \label{Hdelta}%
\end{equation}
and consider the Hamilton-Jacobi-Bellman equation
\begin{equation}
\left\{
\begin{array}
[c]{l}%
\displaystyle\frac{\partial}{\partial t}V^{\delta}\left(  t,x\right)
+\underset{v\in\mathbb{U}}{\inf}H^{\delta}\left(  x,(V^{\delta},\nabla
_{x}V^{\delta},\nabla_{xx}V^{\delta})(t,x),v\right)  =0,\ \left(  t,x\right)
\in\left[  0,T\right]  \times\mathbb{R}^{d},\\
V^{\delta}\left(  T,x\right)  =\Phi_{\delta}(x),\ \ \ x\in\mathbb{R}^{d},
\end{array}
\right.  \label{hjbdelta}%
\end{equation}

Since $H^{\delta}$ is smooth and $(\sigma_{\delta}\sigma_{\delta}^{\ast
})\left(  x,y\right)  $ is uniformly elliptic, then according to
\cite{Krylov1} (Theorems 6.4.3 and 6.4.4), equation \eqref{hjbdelta} admits a
unique solution $V^{\delta}$ which belongs to $C_{b}^{1,2}([0,T]\times
\mathbb{R}^{d})$. The regularity of $V^{\delta}$ and the compactness of the
control state space $\mathbb{U}$ allow us to find a measurable function
$v^{\delta}:[0,T]\times\mathbb{R}^{d}\longmapsto\mathbb{U}$ such that, for all
$(t,x)\in\lbrack0,T]\times\mathbb{R}^{d}$,
\begin{equation}
H^{\delta}\left(  x,(V^{\delta},\nabla_{x}V^{\delta},\nabla_{xx}V^{\delta
})(t,x),v^{\delta}(t,x)\right)  =\underset{v\in\mathbb{U}}{\inf}H^{\delta
}\left(  x,(V^{\delta},\nabla_{x}V^{\delta},\nabla_{xx}V^{\delta
})(t,x),v\right)  . \label{vdelta}%
\end{equation}

Let $(t,x)\in\lbrack0,T]\times\mathbb{R}^{d}$. For $\delta>0$, let $V^{\delta
}$ be the solution of \eqref{hjbdelta} and $v^{\delta}$ the function defined
by \eqref{vdelta}. Consider the SDE:
\begin{equation}
\left\{
\begin{array}
[c]{l}%
dX_{s}^{\delta}\ =\ b_{\delta}(X_{s}^{\delta},V^{\delta}(s,X_{s}^{\delta
}),v^{\delta}(s,X_{s}^{\delta}))ds+\ \sigma_{\delta}(X_{s}^{\delta},V^{\delta
}(s,X_{s}^{\delta}))dW_{s}^{\delta},\,\ s\in\lbrack t,\ T],\\
X_{t}^{\delta}=x.
\end{array}
\right.  \label{sdedelta1}%
\end{equation}
Since $b_{\delta}(x,V^{\delta}(s,x),\nabla_{x}V^{\delta}(s,x)\sigma_{\delta
}(x,V^{\delta}(s,x)),v^{\delta}(s,x))$ and $\sigma_{\delta}(x,V^{\delta
}(s,x))$ are bounded measurable in $(t,x)$ and $\sigma_{\delta}(x,V^{\delta
}(s,x))$ is Lipschitz in $x$ and uniformly elliptic, then according to
\cite{bahlali}, Theorem 2.1 pp 56 (see also \cite{kbcras}), equation
(\ref{sdedelta1}) has a pathwise unique solution $X^{\delta}$.\newline For
$s\in\lbrack t,T]$, let
\begin{equation}
Y_{s}^{\delta}:=V^{\delta}(s,X_{s}^{\delta})\qquad\mbox{  and   }\qquad
Z_{s}^{\delta}:=\nabla_{x}V^{\delta}(s,X_{s}^{\delta})\sigma_{\delta}%
(X_{s}^{\delta},V^{\delta}(s,X_{s}^{\delta})). \label{ydelta=vdelta}%
\end{equation}
Applying Itô's formula to $V^{\delta}(s,X_{s}^{\delta})$, we get:
\begin{equation}
\left\{
\begin{array}
[c]{l}%
dX_{s}^{\delta}=b_{\delta}(X_{s}^{\delta},Y_{s}^{\delta},\ v^{\delta}%
(s,X_{s}^{\delta}))ds+\sigma_{\delta}(X_{s}^{\delta},Y_{s}^{\delta}%
)dW_{s}^{\delta},\quad s\in\lbrack t,\ T]\\
X_{t}^{\delta}=x,\\
dY_{s}^{\delta}=-f_{\delta}(X_{s}^{\delta},Y_{s}^{\delta},Z_{s}^{\delta
},v^{\delta}(s,X_{s}^{\delta}))ds+Z_{s}^{\delta}dW_{s}^{\delta},\\
Y_{T}^{\delta}=\Phi_{\delta}(X_{T}^{\delta})
\end{array}
\right.  \label{ahmed}%
\end{equation}
Since $f_{\delta}$ is uniformly Lipschitz in $(y,z)$, then according to
\cite{PP1}, the backward component of equation (\ref{ahmed}) has a unique
solution $(Y^{\delta},Z^{\delta})$ in $\mathcal{S}_{\nu}^{2}(t,T;\mathbb{R}%
)\times\mathcal{H}_{\nu}^{2}(t,T;\mathbb{R}^{d})$. Therefore $(X^{\delta
},Y^{\delta},Z^{\delta})$ is the unique solution of FBSDE (\ref{ahmed}) in
$\mathcal{S}_{\nu}^{2}(t,T;\mathbb{R}^{d})\times\mathcal{S}_{\nu}%
^{2}(t,T;\mathbb{R})\times\mathcal{H}_{\nu}^{2}(t,T;\mathbb{R}^{d})$.
Therefore, the process $u_{s}^{\delta}:=v^{\delta}(s,X_{s}^{\delta})$ is an
admissible feedback control.

Let $u\in\mathcal{U}_{\nu^{\delta}(t)}$ be an admissible control. Let
$(X^{\delta,t,x,u},Y^{\delta,t,x,u},Z^{\delta,t,x,u})$ be the unique
$\mathcal{F}_{t} $-adapted continuous solution of the following FBSDE defined
on $[t,T]$:
\begin{equation}
\label{xdeltau}\left\{
\begin{array}
[c]{l}%
dX_{s}^{\delta,t,x,u}=b_{\delta}\left(  X_{s}^{\delta,t,x,u},Y_{s}%
^{\delta,t,x,u},u_{s}\right)  ds+\sigma_{\delta}\left(  X_{s}^{\delta
,t,x,u},Y_{s}^{\delta,t,x,u}\right)  dW^{\delta}_{s},\\
X_{t}^{\delta,t,x,u}=x,\\
dY_{s}^{\delta,t,x,u}=-f_{\delta}(X_{s}^{\delta,t,x,u},Y_{s}^{\delta
,t,x,u},Z_{s} ^{\delta,t,x,u},u_{s})ds+Z_{s}^{\delta,t,x,u}dW_{s}^{\delta
}+dM_{s}^{\delta,t,x,u},\\
Y_{T}^{\delta,t,x,u}=\Phi_{\delta}(X_{T}^{\delta,t,x,u}),\\
M^{\delta}\in\mathcal{M}^{2}_{\nu^{\delta}}(t,T;\mathbb{R}^{d})
\mbox{ is orthogonal to } W^{\delta}.
\end{array}
\right.
\end{equation}

The cost functional associated to the controlled FBSDE \eqref{xdeltau} is then
defined by:
\[
J^{\delta}(u):=Y^{\delta,t,x,u}_{t},\, u\in\mathcal{U}_{\nu^{\delta}}(t).
\]

Since $(X^{\delta},Y^{\delta},Z^{\delta})$ satisfies the FBSDE (\ref{xdeltau})
for $u=u^{\delta}$, with $M^{\delta} = 0$, then by the uniqueness of equation
(\ref{xdeltau}), we have $(X^{\delta},Y^{\delta},Z^{\delta})= (X^{\delta
,t,x,u^{\delta}},Y^{\delta,t,x,u^{\delta}},Z^{\delta,t,x,u^{\delta}})$. In
particular $Y^{\delta,t,x,u^{\delta}}_{t}=Y^{\delta}_{t}=V^{\delta}(t,x)$. We
then have proved the following Lemma.

\begin{lemma}
\label{uadmissible} Assume that \textbf{(A)} is satisfied. Then, for any
$\delta> 0$, there exists an admissible feedback control $u^{\delta}_{s}:=
v^{\delta}(s, X^{\delta}_{s})$ defined for $s\in[0,T]$ such that:
\begin{equation}
\label{Jdelta=Vdelta}\displaystyle J^{\delta}(u^{\delta})=V^{\delta}(t,x)=
\hbox{essinf}_{u\in\mathcal{U} _{\nu^{\delta}}(t)}J^{\delta}(u),
\end{equation}

\end{lemma}

The following two lemmas will be needed for the construction of the optimal
control. The second one shows that the variable $Z^{\delta}$ is uniformly
bounded. This allows us to consider $Z^{\delta}$ as a control.

\begin{lemma}
\label{Vdeltacauchy} Assume that $\mathbf{{(A)}}$ is satisfied. Then,

$(i)$ \ there exists a non-negative constant $C$ depending on $K$, $T$ and the
bounds of the coefficients such that,
\begin{equation}
\left\vert V^{\delta}(t,x)-V^{\delta^{\prime}}(t^{\prime},x^{\prime
})\right\vert \leq C(\left\vert \delta-\delta^{\prime}\right\vert
+|x-x^{\prime}|+\left\vert t-t^{\prime}\right\vert ^{1/2}).\label{vari}%
\end{equation}

$(ii)$ \ $V^{\delta}$ converges uniformly to a bounded function $V,$ which is
the unique viscosity solution of the initial HJB equation \eqref{hjb}.
\end{lemma}

\bop$(i)$ \ From the uniqueness of the solution of the controlled forward
equation with control process $u^{\delta},$ it follows that $X^{\delta
,u^{\delta}}=X^{\delta}.$

Let $\delta^{\prime}>0$ and $(t^{\prime},x^{\prime} )\in[0,T]\times
\mathbb{R}^{d}$. For $\delta>0$, let $X^{\delta^{\prime},t^{\prime},x^{\prime
},u^{\delta}}\in\mathcal{S}^{2} _{\nu^{\delta}}(t^{\prime},T;\mathbb{R}^{d})$
denote the unique solution of the forward equation on $[t^{\prime},T]$:
\[
\left\{
\begin{array}
[c]{l}%
dX_{s}^{\delta^{\prime},t^{\prime},x^{\prime},u^{\delta}}= b_{\delta^{\prime}
}\big(X_{s}^{\delta^{\prime},t^{\prime},x^{\prime},u^{\delta}},V^{{\delta
}^{\prime}} (s,X_{s}^{\delta^{\prime},t^{\prime},x^{\prime},u^{\delta}})
,u^{\delta}_{s}\big)ds\\
\hskip 2.3cm + \ \sigma_{\delta^{\prime}}\big( X_{s}^{\delta^{\prime
},t^{\prime},x^{\prime}},V^{{\delta}^{\prime}} (s,X_{s}^{\delta^{\prime
},t^{\prime},x^{\prime},u^{\delta}})\big) dW_{s}^{\delta}\\
X_{t^{\prime}}^{\delta^{\prime},t^{\prime},x^{\prime},u^{\delta}}=x^{\prime}.
\end{array}
\right.
\]
We extend this solution to the whole interval $[0,T]$ by setting
$X_{s}^{\delta^{\prime},t^{\prime},x^{\prime},u^{\delta}}=x^{\prime},$ for
$s<t^{\prime}$. \newline We put for $s\in[t^{\prime},T]$,%

\begin{align}
\label{ftilde}\widetilde{f}^{\delta^{\prime},t^{\prime},x^{\prime},u^{\delta}
}_{s} = -\bigg(  &  \frac{\partial}{\partial s}V^{\delta^{\prime}}%
(s,X^{\delta^{\prime},t^{\prime},x^{\prime},u^{\delta}}_{s})\\
&  + \frac12\mbox{trace} \ \sigma_{\delta^{\prime}}\sigma_{\delta^{\prime}%
}^{\ast} (X^{\delta^{\prime},t^{\prime},x^{\prime},u^{\delta}}_{s}%
,V^{\delta^{\prime}}(s,X^{\delta^{\prime},t^{\prime},x^{\prime},u^{\delta}%
}_{s})) \displaystyle \times\nabla_{xx}V^{\delta} (s,X^{\delta^{\prime
},t^{\prime},x^{\prime},u^{\delta}}_{s})\nonumber\\
&  +b_{\delta^{\prime} }( X_{s}^{\delta^{\prime},t^{\prime},x^{\prime
},u^{\delta}},V^{{\delta}^{\prime}} (s,X_{s}^{\delta^{\prime},t^{\prime
},x^{\prime},u^{\delta}})) \bigg)\nonumber
\end{align}
It{\^o}'s formula applied to \ $V^{\delta^{\prime}}(s,X_{s}^{\delta,t^{\prime
},x^{\prime},u^{\delta}})$ \ shows that the processes
\begin{align*}
&  Y_{s}^{\delta^{\prime},t^{\prime},x^{\prime}}:=V^{\delta^{\prime}}
(s,X_{s}^{\delta^{\prime},t^{\prime},x^{\prime},u^{\delta}}),\\
&  Z_{s} ^{\delta^{\prime},t^{\prime},x^{\prime}}:=\nabla_{x}V^{\delta
^{\prime}}(s,X_{s} ^{\delta^{\prime},t^{\prime},x^{\prime},u^{\delta}})
\sigma_{\delta^{\prime} }(X_{s}^{\delta^{\prime},t^{\prime},x^{\prime
},u^{\delta}},V^{\delta^{\prime}}(s,X_{s} ^{\delta^{\prime},t^{\prime
},x^{\prime},u^{\delta}})),\\
&  M_{s}^{\delta^{\prime}, t^{\prime},x^{\prime}}:=0,\, s\in[t^{\prime},T],
\end{align*}

\noindent is the unique solution of the BSDE
\begin{equation}
\left\{
\begin{array}
[c]{l}%
dY_{s}^{\delta^{\prime},t^{\prime},x^{\prime}}=-\widetilde{f}_{s}%
^{\delta^{\prime},t^{\prime},x^{\prime},u^{\delta}}ds+Z_{s}^{\delta^{\prime
},t^{\prime},x^{\prime}}dW_{s}^{\delta}+dM_{s}^{\delta^{\prime},t^{\prime
},x^{\prime}},\,s\in\lbrack t^{\prime},T],\\
Y_{T}^{\delta^{\prime},t^{\prime},x^{\prime}}=\Phi_{\delta^{\prime}}%
(X_{T}^{\delta^{\prime},t^{\prime},x^{\prime}}),\\
(Y^{\delta^{\prime},t^{\prime},x^{\prime}},Z^{\delta^{\prime},t^{\prime
},x^{\prime}})\in\mathcal{S}_{\nu^{\delta}}^{2}(t^{\prime},T;\mathbb{R}%
)\times\mathcal{H}_{\nu^{\delta}}^{2}(t^{\prime},T;\mathbb{R}^{d}),\\
M^{\delta^{\prime},t^{\prime},x^{\prime}}\in\mathcal{M}_{\nu^{\delta}}%
^{2}(t^{\prime},T;\mathbb{R}^{d})\mbox{ is orthogonal to }W^{\delta}.
\end{array}
\right.  \label{deltap}%
\end{equation}
In order to compute the estimation, let us define the following BSDE%

\begin{equation}
\label{bsdefdel'}\left\{
\begin{array}
[c]{l}%
dY_{s}^{\delta^{\prime},t^{\prime},x^{\prime},u^{\delta}}=-f_{\delta^{\prime}
}\left(  X_{s}^{\delta^{\prime},t^{\prime},x^{\prime},u^{\delta}},Y_{s}
^{\delta^{\prime},t^{\prime},x^{\prime},u^{\delta}}, Z_{s}^{\delta^{\prime
},t^{\prime},x^{\prime},u^{\delta}},u^{\delta}_{s}\right)  ds\\
\qquad\qquad\qquad+Z_{s}^{\delta^{\prime},t^{\prime},x^{\prime},u^{\delta}
}dW^{\delta}_{s} +dM_{s}^{\delta^{\prime},t^{\prime},x^{\prime},u^{\delta}},\,
s\in[t^{\prime},T],\\
\hskip 2mm Y_{T}^{\delta^{\prime},t^{\prime},x^{\prime},u^{\delta}}
=\Phi_{\delta}(X_{T}^{\delta,t^{\prime},x^{\prime},u^{\delta}}),\\
(Y^{\delta^{\prime},t^{\prime},x^{\prime},u^{\delta}},Z^{\delta^{\prime
},t^{\prime},x^{\prime},u^{\delta}})\in\mathcal{S}^{2}_{\nu^{\delta}}
(t^{\prime},T;\mathbb{R}) \times\mathcal{H}^{2}_{\nu^{\delta}}(t^{\prime
},T;\mathbb{R}^{d} ),\\
M^{\delta^{\prime},t^{\prime},x^{\prime},u^{\delta}}\in\mathcal{M}^{2}
_{\nu^{\delta}}(t^{\prime},T;\mathbb{R}^{d}) \mbox{ is orthogonal to }
W^{\delta}.
\end{array}
\right.
\end{equation}
Since $V^{\delta^{\prime}}$ is a classical solution to the
Hamilton-Jacobi-Bellman equation it follows that
\[
\widetilde{f}^{\delta^{\prime},t^{\prime},x^{\prime},u^{\delta}}_{s}\le
f_{\delta^{\prime}}\left(  X_{s}^{\delta^{\prime},t^{\prime},x^{\prime
},u^{\delta}},Y_{s}^{\delta^{\prime},t^{\prime},x^{\prime}}, Z_{s}
^{\delta^{\prime},t^{\prime},x^{\prime}},u^{\delta}_{s}\right)  ,\, \quad
s\in[t^{\prime},T].
\]
Hence, the comparison Theorem shows that \ $Y^{\delta^{\prime},t^{\prime
},x^{\prime}}_{s}\le Y_{s}^{\delta^{\prime},t^{\prime},x^{\prime},u^{\delta}
},$ $s\in[t^{\prime},T]$, $P$-a.s.

\noindent Therefore, we have
\[
V^{\delta^{\prime}}(t^{\prime},x^{\prime})-V^{\delta}(t,x) \leq Y_{t^{\prime}%
}^{\delta^{\prime},t^{\prime},x^{\prime},u^{\delta}} -Y^{\delta,t,x,u^{\delta
}}_{t} ,\, \ \mathbb{P}\mbox{-a.s.}
\]
Using a symmetric argument, we deduce that :
\[
|V^{\delta^{\prime}}(t^{\prime},x^{\prime})-V^{\delta}(t,x)| \leq
|Y_{t^{\prime}}^{\delta^{\prime},t^{\prime},x^{\prime},u^{\delta}}
-Y^{\delta,t,x,u^{\delta}}_{t}| ,\, \ \mathbb{P}\mbox{-a.s.}
\]
Since $V^{\delta}$ and $V^{\delta^{\prime}}$ are deterministic, we have
\[
|V^{\delta^{\prime}}(t^{\prime},x^{\prime})-V^{\delta}(t,x)| \leq
\mathbb{E}(|Y_{t^{\prime}}^{\delta^{\prime},t^{\prime},x^{\prime},u^{\delta}}
-Y^{\delta,t,x,u^{\delta}}_{t}|)
\]
Hence, it suffices to estimate $\mathbb{E}(|Y_{t^{\prime}}^{\delta^{\prime
},t^{\prime},x^{\prime},u^{\delta}} -Y^{\delta,t,x,u^{\delta}}_{t}|)$.

\noindent We assume that $t^{\prime}<t$, \ and \ for $s<t^{\prime}$,
\ \ $Y_{s}^{\delta^{\prime},t^{\prime},x^{\prime},u^{\delta}}=Y_{t^{\prime}%
}^{\delta^{\prime}, t^{\prime},x^{\prime},u^{\delta}},\;\;Z_{s}^{\delta
^{\prime},t^{\prime},x^{\prime},u^{\delta}}=0$ and $M_{s}^{\delta^{\prime
},t^{\prime},x^{\prime},u^{\delta}}=0$.

\vskip 0.2cm\noindent We have,
\begin{align*}
\mathbb{E}(\vert Y_{t}^{\delta,t,x,u^{\delta}} - Y_{t^{\prime}}^{\delta
^{\prime},t^{\prime},x^{\prime},u^{\delta}}\vert^{2})  &  \leq C_{1}%
\mathbb{E}(\vert Y_{t}^{\delta^{\prime},t^{\prime},x^{\prime},u^{\delta}}-
Y_{t^{\prime}}^{\delta^{\prime},t^{\prime},x^{\prime},u^{\delta}}\vert^{2}
+\vert Y_{t}^{\delta,t,x,u^{\delta}} - Y_{t}^{\delta^{\prime},t^{\prime
},x^{\prime},u^{\delta}}\vert^{2} )\\
&  \leq C_{1}(|t-t^{\prime}| + \mathbb{E}[\sup_{t\leq s \leq T}\vert
Y_{s}^{\ delta,t,x,u^{\delta}} - Y_{s}^{\delta^{\prime},t^{\prime},x^{\prime
},u^{\delta}}\vert^{2}] )
\end{align*}
where $C_{1}$ is some positive constant which depends on $T$, $K$ and the
bounds of $b$, $\sigma$, $f$ and $\Phi$ but not on $t, x, \delta$.

It remains to show that:
\begin{align}
\label{DeltaY}\mathbb{E}[\sup_{t\leq s \leq T}\vert Y_{s}^{\delta
,t,x,u^{\delta}} - Y_{s}^{\delta^{\prime},t^{\prime},x^{\prime},u^{\delta}%
}\vert^{2}] ) \leq C(|x-x^{\prime2 }+ |\delta-\delta^{\prime2}),
\end{align}
where $C$ is some positive constant which depends on the $T$ $K$ and the
bounds of $b$, $\sigma$, $f$ and $\Phi$ but not on $t, x, \delta$.

In the sequel of the proof, the positive constant $C$ may be change from line
to line. This constant will depend on the $T$ $K$ and the bounds of $b$,
$\sigma$, $f$ and $\Phi$ but not on $t, x, \delta$. To simplify the notations
throughout this proof, we put
\[
X^{\delta, t,x,u^{\delta}}_{s}:=X^{\delta}_{s},\qquad Y^{\delta,
t,x,u^{\delta}}_{s}:=Y^{\delta}_{s},\qquad Z^{\delta, t,x,u^{\delta}}%
_{s}:=Z^{\delta, }_{s}
\]
and
\[
X^{\delta^{\prime}, t^{\prime},x^{\prime\delta}}_{s} := X^{\delta
^{\prime\delta}}_{s}, \qquad Y^{\delta^{\prime}, t^{\prime},x^{\prime\delta}%
}_{s}:=Y^{\delta^{\prime\delta}}_{s} , \qquad Z^{\delta^{\prime}, t^{\prime
},x^{\prime\delta}}_{s}:=Z^{\delta^{\prime\delta}}_{s}
\]

Since all the coefficients and the terminal data are bounded, then standard
arguments of BSDEs shows that there exists a constant $C>0$ which depends from
$T$ the bounds of $b, \sigma, \Phi$ and $f$ such that
\begin{align}
\label{bornesyyprime}\mathbb{E} \left(  \sup_{t\leq s\leq T}\left[
|X^{\delta}_{s}|^{2} + |X^{\delta^{\prime\delta}}_{s}|^{2} + |Y^{\delta}%
_{s}|^{2} + |Y^{\delta^{\prime\delta}}_{s}|^{2}\right]  + \int_{t}^{T}
|Z^{\delta}_{r}|^{2} dr + \int_{t}^{T} |Z^{\delta^{\prime\delta}}_{r}|^{2}
dr\right)  \leq C
\end{align}
Using Proposition \ref{mollif}, we have
\begin{align}
\label{termin}\mathbb{E}( |Y^{\delta}_{T}- Y^{\delta^{\prime\delta}}_{T}%
|^{2})  &  \leq2\ K^{2} [|\delta-\delta^{\prime2}+ \mathbb{E} ( |X^{\delta
}_{T}-X^{\delta^{\prime\delta}}_{T}|^{2})]\nonumber\\
&  \leq2\ K^{2} [|\delta-\delta^{\prime2}+ \mathbb{E} (\sup_{t\leq s\leq T}
|X^{\delta}_{s}-X^{\delta^{\prime\delta}}_{s}|^{2})]
\end{align}
Using again Proposition \ref{mollif} and standard arguments of BSDEs
(Burkhölder-Davis-Gundy's inequality, etc.) and inequality \eqref{termin} we
show that
\begin{align}
\label{DeltaY2}\mathbb{E} (\sup_{t\leq s\leq T} |Y^{\delta}_{s}-Y^{\delta
^{\prime\delta}} _{s}|^{2})  &  \leq C \bigg(|\delta-\delta^{\prime2}+
\mathbb{E} (\sup_{t\leq s\leq T} |X^{\delta}_{s}-X^{\delta^{\prime\delta}}
_{s}|^{2})\nonumber\\
&  \qquad\qquad+ \int_{t}^{T} \mathbb{E} (\sup_{t\leq s\leq r} |Y^{\delta}%
_{s}-Y^{\delta^{\prime\delta}} _{s}|^{2})dr \bigg)
\end{align}
and
\begin{align}
\label{DeltaX2}\mathbb{E} (\sup_{t\leq s\leq T} |X^{\delta}_{s}-X^{\delta
^{\prime\delta}} _{s}|^{2}) \leq C \bigg(|\delta-\delta^{\prime2}+
|x-x^{\prime2 }+ \int_{t}^{T} \mathbb{E} (\sup_{t\leq s\leq r} |Y^{\delta}%
_{s}-Y^{\delta^{\prime\delta}} _{s}|^{2})dr \bigg)
\end{align}
Inequality \eqref{DeltaY} follows now from inequalities \eqref{bornesyyprime},
\eqref{DeltaY2}, \eqref{DeltaX2} and Gronwall's Lemma. Assertion $(i)$ is proved.






We prove assertion $(ii)$. According to assertion $(i)$, $(V^{\delta})$ is
Cauchy sequence with respect to the uniform convregence norm, in
$(t,x)\in\lbrack0,T]\times\mathbb{R}^{d}$. It then converges uniformly to a
function $\bar{V}$ as $\delta\rightarrow0$. Moreover, inequality
\eqref{bornesyyprime} shows that $V^{\delta}$ is uniformly bounded in
$(t,x,\delta)$. Hence $\bar{V}\in C_{b}([0,T]\times\mathbb{R}^{d})$. Since
$H^{\delta}$ converges uniformly on compact sets to $H$, then using the
stability of viscosity solutions, we get that $\bar{V}$ is a viscosity
solution of equation \eqref{hjb}. Thanks to the uniqueness of the solution of
equation (\ref{hjb}), within the class of continuous function, with at most
polynomial growth, we get that $\bar{V}=V$. This shows that the sequence
$(V^{\delta^{\prime\delta}})$ converges to $V$, as $\delta^{\prime}%
\rightarrow0$. Using inequality \eqref{vari}, we deduce that \ $|V^{\delta
}(t,x)-V(t,x)|\leq C\delta$, \ \ \ for each \ $\delta\in
(0,1]\mbox{ and }(t,x)\in\lbrack0,T]\times\mathbb{R}^{d}.$ \eop


\subsection{The passing to the limits}

We will prove the convergence of the approximating control problem to the
original one. We adapt the idea of \cite{BLRT} to or situation. Put $w_{s}%
^{n}:=\nabla_{x}V^{\delta_{n}}(s,X_{s}^{\delta_{n}})$ and $Z_{s}^{\delta_{n}%
}:=w_{s}^{\delta_{n}}\sigma\left(  X_{s}^{\delta_{n}},Y_{s}^{\delta_{n}%
}\right)  $. Consider the sequence of approximating stochastic controlled
systems $(X^{\delta_{n}},Y^{\delta_{n}},Z^{\delta_{n}},u^{\delta_{n}})$. Since
$u^{\delta_{n}}$ and $w^{\delta_{n}}$ are uniformly bounded, we consider the
couple $(u^{\delta_{n}},w^{\delta_{n}})$ as a relaxed control. We show that
the system $(X^{\delta_{n}},Y^{\delta_{n}},Z^{\delta_{n}},u^{\delta_{n}})$ has
a subsequence which converges in law to some controlled system. And, since we
have assumption \textbf{{(H)}}, we use the result of \cite{enj} to prove that
the limiting process is a strict control.

\begin{proposition}
\label{propbar} Assume that \textbf{{(A)}} and \textbf{{(H)}} are satisfied.
Let $(t,x)\in[0,T]\times\mathbb{R}^{d}$ and $(\delta_{n})_{n\in\mathbb{N}}$ be
a sequence of positive real numbers which tends to $0$. Then, there exists a
reference stochastic system $\bar\nu=(\bar\Omega,\bar{\mathcal{F} }
,\bar{\mathbb{P}},\bar{\mathcal{F}_{t} },\bar W)$, a process $(\bar X,\bar
Y,\bar Z,\bar M)\in\mathcal{S}^{2}_{\bar\nu}(t,T;\mathbb{R}^{d})\times
\mathcal{S}^{2} _{\bar\nu}(t,T;\mathbb{R})\times\mathcal{S}^{2}_{\bar\nu
}(t,T;\mathbb{R} ^{d})\times\mathcal{M}^{2}_{\bar\nu}(t,T;\mathbb{R}^{d})$,
with $\bar M$ orthogonal to $\bar W$, and an admissible control $\bar
u\in\mathcal{U}_{\bar\nu}(t)$, such that: \newline1) There is a subsequence of
$(X^{\delta_{n}},Y^{\delta_{n} })_{n\in\mathbb{N}}$ which converges in
distribution to $(\bar X,\bar Y)$, \newline2) $(\bar X,\bar Y,\bar Z,\bar M)$
is a solution of the following system
\begin{equation}
\label{bsde}\left\{
\begin{array}
[c]{l}%
d\bar X_{s}=b(\bar X_{s},\bar Y_{s},\bar u_{s})ds+\sigma( \bar X_{s},\bar
Y_{s})d\bar W_{s},\\
d\bar Y_{s}=-f(\bar X_{s},\bar Y_{s},\bar Z_{s},\bar u_{s})ds+\bar Z_{s}d\bar
W_{s}+d\bar M_{s},\;\; s\in[t,T]\\
\bar X_{t}=x,\; \bar Y_{T}=\Phi(X_{T}),
\end{array}
\right.
\end{equation}
3) For every $(t,x)\in[0,T]\times\mathbb{R}^{d}$, it holds that
\[
\bar Y_{t}=V(t,x)=\mbox{essinf}_{u\in\mathcal{U}_{\bar\nu}(t)}J(t,x,u)\; ,
\]
i.e. the admissible control $\bar u\in\mathcal{U}_{\bar\nu}(t)$ is optimal for
(\ref{bsde}).
\end{proposition}

The idea of the proof of this theorem consists in the introduction of an
auxiliary sequence of processes (denoted by $(X^{n},Y^{n})$) which satisfies a
forward-system, for each $n$, and for which a relaxed control exists according
to \cite{enj}. We then show that $(X^{n},Y^{n})$ admits a subsequence, which
converges in law to a couple $(\bar{X},\bar{Y})$. Using the convexity
assumption $\mathbf{{(H)}}$, we prove that $(\bar{X},\bar{Y})$ is associated
to a strict control which is optimal for the original control problem. We
finally show that the initial sequence $(X^{\delta_{n}},Y^{\delta_{n}}%
)_{n\in\mathbb{N}}$ and the auxiliary one have the same limits. More
precisely, we define the sequence of auxiliary processes $(X_{s}^{n},Y_{s}%
^{n})$ as the pathwise unique solution of the following controlled forward
system:
\begin{equation}
\left\{
\begin{array}
[c]{l}%
dX_{s}^{n}=b(X_{s}^{n},Y_{s}^{n},u_{s}^{\delta_{n}})ds+\sigma(X_{s}^{n}%
,Y_{s}^{n})dW_{s}^{\delta_{n}},\\
X_{t}^{n}=x,\\
dY_{s}^{n}=-f(X_{s}^{n},Y_{s}^{n},w_{s}^{\delta_{n}}\sigma(X_{s}^{n},Y_{s}%
^{n}),u_{s}^{\delta_{n}})ds+w_{s}^{n}\sigma(X_{s}^{n},Y_{s}^{n})dW_{s}%
^{\delta_{n}}.\\
Y_{t}^{n}=V^{\delta_{n}}(t,x),\ \ s\in\lbrack t,T].
\end{array}
\right.  \label{neqqqq}%
\end{equation}
where $u_{s}^{\delta_{n}}:=v^{\delta_{n}}(s,X_{s}^{\delta_{n}})$ \ and
\ $w_{s}^{\delta_{n}}=\nabla_{x}V^{\delta_{n}}(s,X_{s}^{\delta_{n}})$.


\begin{lemma}
\label{xnxdelta-n} There exists a constant $L>0$ such that for any
$n\in\mathbb{N}$,
\begin{equation}
\label{dddn}%
\begin{array}
[c]{c}%
\mathbb{E}[\sup_{s\in[t,T]}|X_{s}^{\delta_{n}}-X^{n}_{s}|^{2}]\leq L\delta
_{n}^{2},\\
\mathbb{E}[\sup_{s\in[t,T]}|Y_{s}^{\delta_{n}}-Y^{n}_{s}|^{2}]\leq L\delta
_{n}^{2}.
\end{array}
\end{equation}
where $(X^{\delta_{n}}, Y^{\delta_{n}})$ satisfies the FBSDE \eqref{ahmed}.
\end{lemma}


\bop It can be performed as that of Lemma \ref{Vdeltacauchy}, assertion $(i)$.
\eop


\vskip0.2cm\noindent\textbf{Proof of Proposition \ref{propbar}}. Note that for
every $n$, the process $(X_{s}^{\delta_{n}},Y_{s}^{\delta_{n}})$ is a weak
solution to the following controlled forward system:
\begin{equation}
\left\{
\begin{array}
[c]{l}%
dX_{s}^{\delta_{n}}=b_{\delta_{n}}(X_{s}^{\delta_{n}},Y_{s}^{\delta_{n}}%
,u_{s}^{\delta_{n}})ds+\sigma_{\delta_{n}}(X_{s}^{\delta_{n}},Y_{s}%
^{\delta_{n}})dW_{s}^{\delta_{n}},\ \ s\in\lbrack t,T]\\
X_{t}^{\delta_{n}}=x,\\
dY_{s}^{\delta_{n}}=-f_{\delta_{n}}(X_{s}^{\delta_{n}},Y_{s}^{\delta_{n}%
},w_{s}^{n}\sigma_{\delta_{n}}(X_{s}^{\delta_{n}},Y_{s}^{\delta_{n}}%
),u_{s}^{\delta_{n}})ds+w_{s}^{n}\sigma_{\delta_{n}}(X_{s}^{\delta_{n}}%
,Y_{s}^{\delta_{n}})dW_{s}^{\delta_{n}}.\\
Y_{t}^{\delta_{n}}=V^{\delta_{n}}(t,x).
\end{array}
\right.  \label{subneqq}%
\end{equation}
Thanks to \eqref{ydelta=vdelta} we have for $t\leq s\leq T$,
\[
Y_{s}^{\delta_{n}}=V^{\delta_{n}}(s,X^{\delta_{n}}%
)~~~~\mbox{  and     }\ \ \ u_{s}^{\delta_{n}}=v^{\delta_{n}}(s,X_{s}%
^{\delta_{n}}).
\]
Since $(s,x)\mapsto V^{\delta_{n}}(s,x)$ is of class $\mathcal{C}^{1,2}$ and
satisfies equation (\ref{hjb}), then using Itô's formula we get for $t\leq
s\leq T$
\begin{align}
Y_{s}^{\delta_{n}}  &  =\Phi_{\delta_{n}}(X_{T}^{\delta_{n}})+\int_{s}%
^{T}f_{\delta_{n}}(X_{r}^{\delta_{n}},Y_{r}^{\delta_{n}},w_{r}^{\delta_{n}%
}\sigma_{\delta_{n}}(X_{r}^{\delta_{n}},Y_{r}^{\delta_{n}}),u_{r}^{\delta_{n}%
})dr\nonumber\label{itovdelta-n}\\
&  \qquad-\int_{s}^{T}w_{r}^{\delta_{n}}\sigma_{\delta_{n}}(X_{r}^{\delta_{n}%
},Y_{r}^{\delta_{n}})dW_{r}^{\delta_{n}}.
\end{align}
Let $B^{\delta_{n}}$ be a Brownian motion which is independent from
$W^{\delta_{n}}$. If we put
\[
\chi_{s}^{n}:=\left(
\begin{array}
[c]{c}%
X_{s}^{n}\\
Y_{s}^{n}%
\end{array}
\right)  ,\,\ \ \ r_{s}^{n}:=(w_{s}^{n},0,u_{s}^{\delta_{n}}%
)\,\ \ \ \mbox{ and }\,\ \ \ \mathcal{W}^{n}:=\left(
\begin{array}
[c]{c}%
W^{\delta_{n}}\\
B^{\delta_{n}}%
\end{array}
\right)  ,
\]

then the system (\ref{neqqqq}) can be written in the form:
\begin{equation}
\left\{
\begin{array}
[c]{l}%
d\chi_{s}^{n}=\beta(\chi_{s}^{n},r_{s}^{n})ds+\Sigma(\chi_{s}^{n},r_{s}%
^{n})d\mathcal{W}_{s}^{n},\;\ \ \ s\in\lbrack t,T],\\
\chi_{t}^{n}=\left(
\begin{array}
[c]{c}%
x\\
V^{\delta_{n}}(t,x)
\end{array}
\right)  .
\end{array}
\right.  \label{nneq}%
\end{equation}
where $\beta$ and $\Sigma$ are the functions defined in Lemma \ref{H1}.

According to Lemma \ref{Vdeltacauchy}, $w^{n}_{s}:=\nabla_{x}V^{\delta_{n}%
}(s,X^{\delta_{n}}_{s})$ is uniformly bounded. Hence we can interpret
$(r^{n}_{s},s\in[t,T])$ as a control with values in the compact set
\ $A:=\mathbb{U}\times\bar{B}_{C}(0)\times[0,K]$.

In order to pass to the limit in $n$, we inject the controls \ $r^{n}$ \ in
the set of relaxed controls, that is: we consider $r^{n}$ as a random variable
with values in the space $V$ of all Borel measures $q$ on $[0,T]\times A$,
whose projection $q(\cdot\times A)$ concides with the Lebesgue measure. To
this end, we identify the control process \ $r^{n}$ \ with the random measure
\begin{equation}
q^{n}(\omega,ds,da)=\delta_{r_{s}^{n}(\omega)}(da)ds,\;(s,a)\in\lbrack
0,T]\times A,\omega\in\Omega. \label{qn}%
\end{equation}
From the boundedness of $\{\left(  \Sigma(x,y,z,\theta),\beta(x,y,z,\theta
,v)\right)  ,(x,y,z,\theta,v)\in\mathbb{R}^{d}\times\mathbb{R}\times A\}$ and
the compactness of $\mathcal{V}$, with respect to the topology induced by the
weak convergence of measures, we get the tightness of the laws of $(\chi
^{n},q^{n}),n\geq1,$ on $C([0,T];\mathbb{R}^{d}\times\mathbb{R})\times
\mathcal{V}$. Therefore we can find a probability measure $Q$ on
$C([0,T];\mathbb{R}^{d}\times\mathbb{R})\times\mathcal{V}$ and extract a
subsequence, still denoted by $(\chi^{n},q^{n})$, which converges in law to
the canonical process $(\chi,q)$ on the space $C([0,T];\mathbb{R}^{d}%
\times\mathbb{R})\times\mathcal{V}$ endowed with the measure $Q$. \newline
Since the coefficients of the system (\ref{nneq}) satisfy assumption
$\mathbf{{(H)}}$, then, according to \cite{enj}, there exists a stochastic
reference system $\bar{\nu}=(\bar{\Omega},\bar{\mathcal{F}_{t}},\bar
{\mathbb{P}},\bar{\mathcal{F}_{t}},\bar{\mathcal{W}})$ enlarging
$(C([0,T];\mathbb{R}^{d}\times\mathbb{R})\times V;Q)$ and an $\bar
{\mathcal{F}_{t}}$-adapted process $(\chi,\bar{r})$ ($\bar{r}$ with values in
$A$), which satisfies
\begin{equation}
\left\{
\begin{array}
[c]{l}%
d\chi_{s}=\beta(\chi_{s},\bar{r}_{s})ds+\Sigma(\chi_{s},\bar{r}_{s}%
)d\bar{\mathcal{W}}_{s},\;s\in\lbrack t,T],\\
\chi_{t}=\left(
\begin{array}
[c]{c}%
x\\
V(t,x)
\end{array}
\right)  .
\end{array}
\right.  \label{chi}%
\end{equation}
Moreover, $\chi$ has the same law under $\bar{\mathbb{P}}$ as under $Q$.

If we set $\chi:=\left(
\begin{array}
[c]{c}%
\bar{X}\\
\bar{Y}%
\end{array}
\right)  $, $\bar{\mathcal{W}}:=\left(
\begin{array}
[c]{c}%
\bar{W}\\
\bar{B}%
\end{array}
\right)  $ and $\bar{r}:=(\bar{w},\bar{\theta},\bar{u})$, then the system
\eqref{chi} can be written as follows:
\[
\left\{
\begin{array}
[c]{l}%
d\bar{X}_{s}=b(\bar{X}_{s},\bar{Y}_{s},\bar{u}_{s})ds+\sigma(\bar{X}_{s}%
,\bar{Y}_{s})d\bar{W}_{s},\\
dY_{s}=-f(\bar{X}_{s},\bar{Y}_{s},\bar{Z}_{s},\bar{u}_{s})ds+\bar{Z}_{s}%
d\bar{W}_{s}+\bar{\theta}_{s}d\bar{B}_{s},\;\;s\in\lbrack t,T]\\
\bar{X}_{t}=x,\;\bar{Y}_{t}=V(t,x).
\end{array}
\right.
\]
This proves assertion 1.

Let us prove prove assertion 2). Lemma \ref{xnxdelta-n} shows that if the
sequence $(X^{n},Y^{n})_{n\in\mathbb{N}}$ converges in law, then the same
holds true for $(X^{\delta_{n}},Y^{\delta_{n}})_{n\in\mathbb{N}}$, and the
limits have the same law. Further, we deduce from (\ref{dddn}) and Lemma
\ref{Vdeltacauchy}, that $\bar{Y}_{s}=V(s,\bar{X}_{s}),$ for each $s\in\lbrack
t,T],$ $\bar{\mathbb{P}}-$a.s. In particular, $Y_{T}=\Phi(X_{T}),$
$\bar{\mathbb{P}}$-a.s. Thus, if we set $\bar{M}_{s}=\int_{t}^{s}\bar{\theta
}_{r}d\bar{B}_{r}$, then $\langle\bar{M},\bar{W}\rangle_{s}=\int_{t}^{s}%
\bar{\theta}_{r}d\langle\bar{B},\bar{W}\rangle_{r}=0$ and $(\bar{X},\bar
{Y},\bar{Z},\bar{M})$ satisfies (\ref{bsde}). Assertion 2) is proved.

Let us prove assertion 3). We have already seen that $\bar{Y}_{s}=V(s,\bar
{X}_{s})$ for all $s\in\lbrack t,T],$ $\bar{\mathbb{P}}$-a.s. On the other
hand, according to Li and Wei \cite{LW}, the unique bounded viscosity solution
$V$ of the Hamilton-Jacobi-Bellman equation (\ref{hjb}) satisfies,
\[
V(t,x)=\mbox{essinf}_{u\in\mathcal{U}_{\bar{\nu}^{\delta}}(t)}J(t,x,u),\;\bar
{\mathbb{P}}\mbox{-a.s.}
\]
Assertion 3 is proved. \eop

\begin{remark}
$(i)$ \ As explained in introduction, the uniform Lipschitz condition is not
sufficient to guarantee the existence of solutions and hence the existence of
optimal controls fails also.

$(ii)$ \ When the coefficients $\sigma$ and/or $b$ depend also from the
$z$--variable, the existence and uniqueness of solutions has been established
in the \cite{del} for FBSDEs with non degenerate diffusion and uniformly
Lipschitz coefficients. But the existence of an optimal control, in this case,
seems difficult to obtain. However if we replace the non degeneracy condition
on $\sigma$ by the so called $G$--monotony condition on the coefficients
introduced in \cite{PW}, the existence of an optimal control can be obtained
even when the coefficient $b$ depend from the $z$--variable and the control
$u$ enter the diffusion coefficient $\sigma$. This is the goal of the
forthcoming paper \cite{bkm}.

$(iii)$ \ There are some recent results on the existence and uniqueness of
solutions to fully coupled FBSDEs where all the coefficients depend from the
three variables $x$, $y$ and $z$ under the uniform Lipschitz condition and
supplementary assumptions on the coefficients, see \cite{maal, WY, Z}. But in
the best of our knowledge, the existence of an optimal control under the
assumptions of \cite{maal, WY, Z} is not known.

In our opinion the existence of an optimal control under the conditions used
in \cite{Z} can be obtained by using the method we develop here.

One challenging problem consists to establish the existence of an optimal
control for a fully FBSDE when the coefficient $\sigma$ depends from $z$ and
$u$. In this case, the existence of solutions follows from \cite{WY} and the
Bellman dynamic programming principle is given \cite{LW}.
\end{remark}


\section{Appendix}

\subsection*{ Appendix A : Convexity hypothesis}

\textbf{Proof of Lemma \ref{H1}.} \ Let $\mu$ be a probability measure on the
set $\mathbb{U}\times\bar{B}_{C}(0)$. Our goal is to find a triplet $(\bar
{w},\bar{\theta},\bar{u})\in\mathbb{R}^{d}\times\lbrack0,K]\times\mathbb{U}$
which satisfies :
\begin{equation}%
\begin{array}
[c]{r}%
\int_{\mathbb{U}\times\bar{B}_{C}(0)}((\Sigma\Sigma^{\ast})(x,y,w,0),\beta
(x,y,w)\mu(du,dw)\\
=\left(  (\Sigma\Sigma^{\ast})(x,y,\bar{w},\bar{\theta}),\beta(x,y,\bar
{w},\bar{\theta},\bar{u})\right)  .
\end{array}
\label{tripp}%
\end{equation}
Let $\Phi(u,w)=\left(  (\sigma\sigma^{\ast})(x,y),w\sigma\sigma^{\ast
}(x,y),b(x,y,u),f(x,y,w\sigma(x,y),u)\right)  $. According to assumption
$\mathbf{{(H)}}$ and the continuity of $\Phi,$ there exists $(\bar{u},\bar
{w})$ in $\mathbb{U}\times\bar{B}_{C}(0)$ such that
\begin{equation}
\int_{\mathbb{U}\times\bar{B}_{C}(0)}\Phi(u,w)\mu(du,dw)=\Phi(\bar{u},\bar
{w}). \label{vww}%
\end{equation}
A simple computation gives,
\[
\Sigma\Sigma^{\ast}(x,y,w,\theta)=%
\begin{pmatrix}
\sigma\sigma^{\ast}(x,y) & \sigma\sigma^{\ast}(x,y)w^{\ast}\\
w\sigma\sigma^{\ast}(x,y) & w\sigma\sigma^{\ast}(x,y)w^{\ast}+\theta^{2}%
\end{pmatrix}
\]
The expression of $(\Sigma\Sigma^{\ast})(x,y,w,0)$ shows that, to obtain
(\ref{tripp}), it suffices to find $\bar{\theta}\in\lbrack0,K]$ such that
\begin{equation}
\bar{\theta}^{2}=\int_{\mathbb{U}\times\bar{B}_{C}(0)}w\sigma\sigma^{\ast
}(x,y)w^{\ast}\mu(du,dw)-\bar{w}\sigma\sigma^{\ast}(x,y)\bar{w}^{\ast}%
:=\alpha. \label{pii}%
\end{equation}
Since $\sigma\sigma^{\ast}(x,y,\bar{u})=\int_{\mathbb{U}\times\bar{B}_{C}%
(0)}\sigma\sigma^{\ast}(x,y)\mu(du,dw)$, then we can write $\alpha$ as
follows
\begin{align}
\alpha &  =\int_{\mathbb{U}\times\bar{B}_{C}(0)}w\sigma\sigma^{\ast
}(x,y)w^{\ast}\mu(du,dw)-\int_{\mathbb{U}\times\bar{B}_{C}(0)}\bar{w}%
\sigma\sigma^{\ast}(x,y)\mu(du,dw)\bar{w}^{\ast}\\
&  =\int_{\mathbb{U}\times\bar{B}_{C}(0)}((w-\bar{w})\sigma(x,y))((w-\bar
{w})\sigma(x,y))^{\ast}\mu(du,dw)
\end{align}
It follows that $\alpha\geq0.$ Hence, it suffices now to choose $\bar{\theta
}=\sqrt{\alpha}$.

\noindent Now, from (\ref{pii}) we have
\[
\int_{\mathbb{U}\times\bar{B}_{C}(0)}|w\sigma(x,y)|^{2}\mu(du,dw)=|\bar
{w}\sigma(x,y)|^{2}+\bar{\theta}^{2}.
\]
Since $|\sigma(x,y)|$ is bounded and the support of $\mu$ is included in
$\mathbb{U}\times\bar{B}_{C}(0)$, it follows that $\bar{\theta}$ is bounded,
that is: there exists $K>0$ such that $\bar{\theta}$ belongs to $[0,K]$.
\eop

\end{document}